\documentclass[12pt,a4paper,reqno]{amsart}
\newcommand{\guio}[1]{\nobreakdash-\hspace{0pt}#1}

\usepackage{diagrams,amssymb}
\diagramstyle[height=2em,width=2em, PostScript]
\newarrow{Epi}----{>>}
\newarrow{Mono}>--->
\newarrow{Iso}>---{>>}
\newarrow{Inc}C--->
\newarrow{Mapsto}|--->
\newarrow{Igual}=====
\newarrow{Dashto}{}{dash}{}{dash}>
\newcommand{\Z}{{\mathbb Z}}

\newcommand{\C}{{\mathbb C}}
\newcommand{\F}{{\mathbb F}}

\newcommand{\map}{\operatorname{map}\nolimits}

                     \newcommand{\Rep}{\operatorname{Rep}\nolimits}
\newcommand{\A}{\ifmmode{\mathcal{A}}\else${\mathcal{A}}$\fi}
\newcommand{\K}{\ifmmode{\mathcal{K}}\else${\mathcal{K}}$\fi}
\newcommand{\U}{\ifmmode{\mathcal{U}}\else${\mathcal{U}}$\fi}
\newcommand{\T}{\ifmmode{\mathcal{T}}\else${\mathcal{T}}$\fi}
\newcommand{\FF}{\ifmmode{\mathcal{F}}\else${\mathcal{F}}$\fi}
\newcommand{\LL}{\ifmmode{\mathcal{L}}\else${\mathcal{L}}$\fi}
\newtheorem{Thm}{Theorem}[section]
\newtheorem{Prop}[Thm]{Proposition}
\newtheorem{Cor}[Thm]{Corollary}
\newtheorem{Lem}[Thm]{Lemma}

\theoremstyle{definition}

\newtheorem{Ex}[Thm]{Example}

\theoremstyle{remark}

\title[On the homotopy groups of $p$-completed classifying spaces]
{On the homotopy groups of \boldmath$p$-completed classifying spaces}

\author{Nat\`{a}lia Castellana}
\author{Juan A. Crespo}
\author{J\'er\^{o}me Scherer}

\thanks{All three authors are partially supported by MEC grant MTM2004-06686.
The third author is supported by the program Ram\'on y Cajal, MEC,
Spain, and thanks the CIB (Centre Interfacultaire Bernoulli),
EPFL, Lausanne for its hospitality.}

\begin{document}

%%%%%%%%%%%%%%%%%%%%%%%%%%%%%%%%%%%%%%

\begin{abstract}
Among the generalizations of Serre's theorem on the homotopy
groups of a finite complex we isolate the one proposed by Dwyer
and Wilkerson. Even though the spaces they consider must be
$2$-connected, we show that it can be used to both recover known
results and obtain new theorems about $p$-completed classifying
spaces.
\end{abstract}

%%%%%%%%%%%%%%%%%%%%%%%%%%%%%%%%%%%%%%

\maketitle

%%%%%%%%%%%%%%%%%%%%%%%%%%%%%%%%%%%%%%

%%%%%%%%%%%%%%%%%%%%%%%%%%%%%%%%%%%%%%%%%%%%%%%%%
%%%%%%%%%%%%%%%%%%%%%%%%%%%%%%%%%%%%%%%%%%%%%%%%%
\section*{Introduction}
\label{sec intro}
%%%%%%%%%%%%%%%%%%%%%%%%%%%%%%%%%%%%%%%%%%%%%%%%%
%%%%%%%%%%%%%%%%%%%%%%%%%%%%%%%%%%%%%%%%%%%%%%%%%
In 1953 Serre proved in his celebrated paper~\cite{MR0060234} that
a simply connected finite $CW$\guio{complex} has infinitely many
non-trivial homotopy groups. He conjectured that it should
actually have infinitely many non-trivial homotopy groups with
$2$-torsion, which was proved by McGibbon and Neisendorfer in
1983~\cite{MR749108} by using Miller's solution~\cite{Miller} of
the Sullivan conjecture. They show this phenomenon holds for any
simply connected $CW$-complex with finite mod $2$ cohomology,
replacing thereby the geometric finiteness condition by a purely
algebraic one. Later, in 1986, Lannes and Schwartz~\cite{MR827370}
were able to relax the finiteness condition to \emph{locally
finite} mod $p$ cohomology, i.e. the cohomology is a direct limit
of finite unstable modules over the Steenrod algebra. So they
proved Serre's conjecture for ``Miller spaces", that is,
$1$-connected spaces $X$ for which the space of pointed maps from
$B\Z/p$ to $X$ is contractible.

In 1990, Dwyer and Wilkerson~\cite{MR92b:55004} show the following
generalization of Serre's conjecture. Let $X$ be a $2$-connected
$CW$-complex  of finite type with non-trivial mod $p$ cohomology,
and such that the module of indecomposable elements in
$H^*(X;\F_p)$ is locally finite. Then, infinitely many homotopy
groups of $X$ contain $p$-torsion. This \emph{new} algebraic
condition obviously includes the previous ones, namely spaces with
finite or locally finite mod $p$ cohomology. Moreover, their
condition enables to study spaces with finitely generated mod $p$
cohomology, because the module of indecomposable elements is then
finite.

\pagebreak

The cost in the Dwyer-Wilkerson theorem is that one has to work
with $2$-connected spaces. As they say, ``the example of $\C
P^\infty$ shows that it would not be enough to assume that $X$ is
$1$-connected". This is basically the only simply connected
Postnikov piece with locally finite module of indecomposable
elements, compare with Grodal's \cite[Theorem~1.2]{MR1622342}.

\medskip
\noindent
%%%%%%%%%%%%%%%%%%%%%%%%%%%%%%%%%%%%%%%%%%%%%%%%%
\textbf{Theorem~\ref{thm unificador}.}
%%%%%%%%%%%%%%%%%%%%%%%%%%%%%%%%%%%%%%%%%%%%%%%%%
%\noindent
{\it Let $X$ be a $p$-complete space such that $H^*(X; \F_p)$ is
of finite type. Assume that the module of indecomposable elements
$QH^*(X;\F_p)$ is locally finite. Then one of the following
properties is satisfied:

\begin{itemize}
\item[(1)] $X$ is aspherical, \item[(2)] $X\langle 1 \rangle$ is a
$K(\Z^\wedge_p, 2)^n$, \item[(3)] $X$ has infinitely many homotopy
groups with $p$-torsion.

\end{itemize}

\noindent In the last case the space $\Omega X$ has infinitely
many non-trivial $k$-invariants. }

\medskip

In particular, this implies the Lannes-Schwartz result,
see~Corollary~\ref{cor LS}, and in fact Theorem~\ref{thm
unificador} can even be applied to understand spaces having a
non-trivial fundamental group, such as classifying spaces of
discrete groups. A very exciting problem in homotopy theory is to
determine the behavior of the $p$-completion of classifying
spaces. When $G$ is a finite group, Levi proves
in~\cite{MR1308466} that either $(BG)^\wedge_p$ is again an
Eilenberg-Mac~Lane space or it has infinitely many non-trivial
homotopy groups. Later, Bastardas and Descheemaker discovered the
same phenomenon holds for any virtually nilpotent group. This is
done in~\cite{MR1942236} for torsion free groups and the general
case is solved in~\cite{tesisgemma}. We show that all these
results can be deduced from Theorem~\ref{thm unificador} and we
obtain the same statement for certain quasi $p$-perfect groups of
finite virtual mod~$p$ cohomology and also for the new concept of
$p$-local finite group, due to Broto, Levi, and
Oliver~\cite{MR1992826}.

\medskip
\noindent
%%%%%%%%%%%%%%%%%%%%%%%%%%%%%%%%%%%%%%%%%%%%%%%%%
\textbf{Theorem~\ref{groupdichotomy}.}
%%%%%%%%%%%%%%%%%%%%%%%%%%%%%%%%%%%%%%%%%%%%%%%%%
%\noindent
{\it Let $X$ be the classifying space of a member of the
 following four families:
\begin{itemize}
\item[(1)]  finite groups,

\item[(2)]  $p$-local finite groups,

\item[(3)]  finitely generated virtually nilpotent groups,

\item[(4)]  quasi $p$-perfect groups of finite virtual mod $p$
cohomology.

\end{itemize}

Then the $p$-completion of $X$ is either aspherical or it has
infinitely many homotopy groups with $p$-torsion. In this case the
space $\Omega (X^\wedge_p)$ has infinitely many non-trivial
$k$\guio{invariants}. }

\medskip

%%%%%%%%%%%%%%%%%%%%%%%%%%%%%%%%%%%%%%%%%%%%%%%%%
%%%%%%%%%%%%%%%%%%%%%%%%%%%%%%%%%%%%%%%%%%%%%%%%%
\section{Local loop spaces}
\label{sec local}
%%%%%%%%%%%%%%%%%%%%%%%%%%%%%%%%%%%%%%%%%%%%%%%%%
%%%%%%%%%%%%%%%%%%%%%%%%%%%%%%%%%%%%%%%%%%%%%%%%%
The grounding result for this paper is the equivalence between the
algebraic condition that the module of indecomposable elements
$QH^* (X; \F_p)$ be locally finite and the topological one that
the loop space $\Omega X$ is $B\Z/p$-local. The proof of
\cite[Theorem~3.2]{MR92b:55004} is done at the prime $2$ for
$1$-connected spaces. Although this is probably well-known to the
experts, we give here an alternative proof for this result that
includes arbitrary connected spaces.

%%%%%%%%%%%%%%%%%%%%%%%%%%%%%%%%%%%%%%%%%%%%%%%%%
\begin{Lem}
\label{lem DWbetter}
%%%%%%%%%%%%%%%%%%%%%%%%%%%%%%%%%%%%%%%%%%%%%%%%%
Let $X$ be a $p$-complete, connected space such that $H^*(X;
\F_p)$ is of finite type. Then $QH^*(X; \F_p)$ is locally finite
if and only if $\Omega X$ is $B\Z/p$-local.
\end{Lem}

\begin{proof}
By \cite[Proposition~3.9.7 and 6.4.5]{MR95d:55017} $QH^*(X; \F_p)$
is locally finite if and only if\linebreak $T_V (H^*(X;
\F_p))_{T_V(c)} \cong H^*(X; \F_p)$ for any elementary abelian
$p$-group $V$. By \cite[Proposition~3.4.4]{MR93j:55019}, this is
so if and only if Lannes' $T$ functor computes the cohomology of
$\map(BV, X)_c.\!$ Therefore, the above isomorphism can be
restated by saying that $\map(BV, X)_c \simeq X$, i.e. $\Omega X$
is $B\Z/p$-local.
\end{proof}

Observe that many interesting spaces verify the condition that the
module of indecomposable elements is locally finite. Let us
mention one class of examples taken from
Lannes,~\cite{MR93j:55019}. Recall that a group verifies
\emph{virtually} a certain property if it admits a subgroup of
finite index which verifies the property.

%%%%%%%%%%%%%%%%%%%%%%%%%%%%%%%%%%%%%%%%%%%%%%%%%
\begin{Prop}\cite[p. 203]{MR93j:55019}
\label{prop Lannes}
%%%%%%%%%%%%%%%%%%%%%%%%%%%%%%%%%%%%%%%%%%%%%%%%%
Let $G$ be a group of virtually finite mod $p$ cohomological
\linebreak dimension. Then there is an isomorphism $T H^* (G;
\F_p) \cong \prod_{\rho \in Rep(\Z/p, G)} H^* (C_G(\rho); \F_p)$.
If moreover $G$ has virtually finite mod $p$ cohomology, then
there is a weak equivalence of mapping spaces $\map(B\Z/p,
BG)^\wedge_p \simeq \map(B\Z/p, (BG)^\wedge_p)$. \hfill{\qed}
\end{Prop}

The second part of the proposition is a consequence of
\cite[Proposition~3.4.3]{MR93j:55019} which states that the
$T$-functor computes the cohomology of the corresponding mapping
space.

%%%%%%%%%%%%%%%%%%%%%%%%%%%%%%%%%%%%%%%%%%%%%%%%
\begin{Cor}
\label{cor Lannes}
%%%%%%%%%%%%%%%%%%%%%%%%%%%%%%%%%%%%%%%%%%%%%%%%
Let $G$ be a group of virtually finite mod $p$ cohomology. Then,
the space $\Omega\bigl( (BG)^\wedge_p \bigr)$ is $B\Z/p$-local.
\end{Cor}

\begin{proof}
Consider the fibration $\map_*(B\Z/p,  (BG)^\wedge_p) \rightarrow
\map(B\Z/p, (BG)^\wedge_p) \rightarrow (BG)^\wedge_p$. By
Proposition~\ref{prop Lannes}, we know that the total space is
equivalent to the $p$-completed mapping space $\map(B\Z/p,
BG)^\wedge_p$, which can be in turn identified with the
$p$-completion of $\coprod BC_G(\rho)$, where the disjoint union
is taken over the representations $\Rep(\Z/p, G)$, see for example
\cite[Proposition~7.1]{MR1924024}. The base point is given by the
trivial representation, i.e. lies in the classifying space of the
trivial representation, whose centralizer is $G$ itself.
Therefore, looping once the above fibration, we obtain that
$\map_*(B\Z/p, \Omega(BG)^\wedge_p)$ is contractible.
\end{proof}

%%%%%%%%%%%%%%%%%%%%%%%%%%%%%%%%%%%%%%%%%%%%%%%%%
\begin{Ex}
\label{ex virtuallynilpotent}
%%%%%%%%%%%%%%%%%%%%%%%%%%%%%%%%%%%%%%%%%%%%%%%%%
A virtually nilpotent group $G$ is by definition an extension of a
finite group $Q$ by a nilpotent group $N$. We notice first that if
$G$ is finitely generated one can always assume that $N$ is
torsion free since any finitely generated nilpotent group is
virtually torsion free. Finitely generated torsion free nilpotent
groups have finite cohomological dimension, see for example
\cite[VIII.2]{MR672956}. We infer from the above corollary that
$\Omega(BG)^\wedge_p$ is $B\Z/p$-local for any finitely generated
virtually nilpotent group $G$.

Consider the inclusion $O^p(Q) \rightarrow Q$ of the
maximal $p$-perfect subgroup of $Q$ as in~\cite{MR1308466}
and construct the pull-back $G' = \lim(O^p(Q) \rightarrow Q \leftarrow G)$.
Since the quotient $Q/O^p(Q)$ is a $p$-group $P$, the fibration
$BO^p(Q) \rightarrow BQ \rightarrow BP$ is preserved by
$p$-completion, and so is the pull-backed one $BG' \rightarrow BG
\rightarrow BP$. Therefore $(BG)^\wedge_p$ is the total space
of a fibration
$$
(BG')^\wedge_p \rTo (BG)^\wedge_p \rTo BP
$$
where $P$ is a finite $p$-group and $G'$ has a normal,
finitely generated, torsion-free, nilpotent subgroup $N$ such that
the quotient is $p$-perfect.

When $G$ is virtually nilpotent, finitely generated, and torsion
free, $BG$ is an {\em infra-nilmanifold} (see
\cite[Theorem~3.1.3]{MR1482520}) so that the cohomology of $G$
itself is finite dimensional. In this case $(BG)^\wedge_p$ is
$B\Z/p$-local (as is its loop space of course).
\end{Ex}

The class of groups of finite virtual cohomological dimension is
much larger than the class of virtually nilpotent ones, but we do
not know if all these groups are $\F_p$-good, which prevents us
from being able to obtain our results in full generality. It is
well-known that spaces with $p$-perfect fundamental group are
$\F_p$-good, \cite[Proposition~3.2]{MR51:1825}.

%%%%%%%%%%%%%%%%%%%%%%%%%%%%%%%%%%%%%%%%%%%%%%%%%
\begin{Ex}
\label{ex finitevcd}
%%%%%%%%%%%%%%%%%%%%%%%%%%%%%%%%%%%%%%%%%%%%%%%%%
Let $G$ be a $p$-perfect group of virtually finite mod $p$
cohomology. Then $(BG)^\wedge_p$ is simply connected and $\Omega
\bigl( (BG)^\wedge_p \bigr)$ is $B\Z/p$-local. Examples of such
groups are given by the special linear groups $SL_n(\Z)$ and the
Steinberg groups $St_n(\Z)$, which are even perfect groups. The
homotopy groups of their $p$-completed classifying spaces are
closely related to the algebraic $K$-theory groups of $\Z$.
\end{Ex}

Our next example is a slight generalization. Recall that the
\emph{lower $p$-central series} of a group $G$ is defined
inductively by $\Gamma_0^p(G)=G$ and $\Gamma_{n+1}^p(G)$ is
generated by elements of form $xyx^{-1}y^{-1}z^p$ for $x \in G$
and $y, z \in \Gamma_n^p(G)$. In particular, a group $G$ is
$p$-perfect if and only if $G = \Gamma_1^p(G)$.

%%%%%%%%%%%%%%%%%%%%%%%%%%%%%%%%%%%%%%%%%%%%%%%%%
\begin{Ex}
\label{ex quasiperfect}
%%%%%%%%%%%%%%%%%%%%%%%%%%%%%%%%%%%%%%%%%%%%%%%%%
In analogy with the terminology used by Wagoner in
\cite{MR0354816} and Loday in \cite{MR0447373}, we say that a
group $G$ is \emph{quasi $p$-perfect} if the subgroup
$\Gamma_1^p(G)$ is $p$-perfect. This means that $G$ is an
extension of an elementary abelian $p$-group with a $p$-perfect
one. To make sure that $BG$ is $\F_p$-good we impose the following
condition:

For any finite set $g_1, \dots , g_n$ of elements in
$\Gamma_1^p(G)$ and $g \in G$ there exists an element~$h\! \in\!
\Gamma_1^p(G)$ such that $g g_i g^{-1} = h g_i h^{-1}$ for all $1
\leq i \leq n$.

This actually turns $(BG)^\wedge_p$ into a simple space, compare
with \cite[Lemma~1.3]{MR0354816}. If one requires that $G$ has
virtually finite mod $p$ cohomology, one obtains new examples of
groups~$G$ such that $\Omega \bigl( (BG)^\wedge_p \bigr)$ is
$B\Z/p$-local.
\end{Ex}

%%%%%%%%%%%%%%%%%%%%%%%%%%%%%%%%%%%%%%%%%%%%%%%%%
\begin{Ex}
\label{ex plocal}
%%%%%%%%%%%%%%%%%%%%%%%%%%%%%%%%%%%%%%%%%%%%%%%%%
Let $(S, \FF, \LL)$ be a $p$-local finite group, as defined by
Broto, Levi, and Oliver in~\cite[Definition~1.8]{MR1992826} and
consider its classifying space $|\LL|^\wedge_p$. We know from
\cite[Theorem~5.8]{MR1992826} that $H^*(|\LL|^\wedge_p; \F_p)$ is
noetherian (it can be computed in fact by stable elements, just
like the cohomology of an ordinary finite group). Therefore, the
module of indecomposable elements is finite, and hence $\Omega
(|\LL|^\wedge_p)$ is $B\Z/p$-local by Lemma~\ref{lem DWbetter}.
\end{Ex}

%%%%%%%%%%%%%%%%%%%%%%%%%%%%%%%%%%%%%%%%%%%%%%%%%
%%%%%%%%%%%%%%%%%%%%%%%%%%%%%%%%%%%%%%%%%%%%%%%%%
\section{The Dwyer-Wilkerson theorem}
\label{sec DW}
%%%%%%%%%%%%%%%%%%%%%%%%%%%%%%%%%%%%%%%%%%%%%%%%%
%%%%%%%%%%%%%%%%%%%%%%%%%%%%%%%%%%%%%%%%%%%%%%%%%

We recall in this section the theorem of Dwyer and Wilkerson about
homotopy groups of $2$-connected spaces with locally finite module
of indecomposable elements (their statement is about CW-complexes,
but it holds under the more general assumptions of
\cite[Theorem~1.2]{MR92b:55004}). We explain then how it can be
efficiently applied to understand certain spaces which are not
$2$-connected by considering their $2$-connected cover.

%%%%%%%%%%%%%%%%%%%%%%%%%%%%%%%%%%%%%%%%%%%%%%%%%
\begin{Thm}\cite[Theorem~1.3]{MR92b:55004}
\label{thm DW}
%%%%%%%%%%%%%%%%%%%%%%%%%%%%%%%%%%%%%%%%%%%%%%%%%
Let $X$ be a $2$-connected space such that the mod $p$ cohomology
$H^*(X; \F_p)$ is of finite type. Assume that $\tilde H^*(X; \F_p)
\neq 0$ and that the module of indecomposable elements
$QH^*(X;\F_p)$ is locally finite. Then there exist infinitely many
integers $k$ such that $\pi_k X$ contains $p$-torsion.
\hfill{\qed}
\end{Thm}

The following elementary lemma (compare with
\cite[Lemma~1.4.4]{MR827370}) is the key to understand which are
the spaces that make it impossible to relax the connectivity
assumption in the theorem.

%%%%%%%%%%%%%%%%%%%%%%%%%%%%%%%%%%%%%%%%%%%%%%%%%
\begin{Lem}
\label{lem KZ2}
%%%%%%%%%%%%%%%%%%%%%%%%%%%%%%%%%%%%%%%%%%%%%%%%%
Let $X=K(A, 2)$ be a $p$-complete space such that $H^*(X; \F_p)$
is of finite type. Assume that $\Omega X$ is $B\Z/p$-local. Then
$A$ is isomorphic to a finite direct sum of copies of
$\Z^\wedge_p$.
\end{Lem}

\begin{proof}
Obviously $A$ must be $p$-torsion free since $\Omega X \simeq K(A,
1)$ is assumed to be $B\Z/p$-local. Since $H^*(X; \F_p)$ is of
finite type, we infer that $H_1(A; \F_p)$, which is isomorphic to
$H_2(X; \F_p)$, is finite. Hence by \cite[Lemma~7.5]{MR1062866}
$A$ is an abelian $p$-torsion free Ext-$p$-complete group of
finite type, i.e. $A$ is isomorphic to a finite direct sum of
copies of $\Z^\wedge_p$ (use Harrison's classification
\cite[VI.4.5]{MR51:1825} or Bousfield's comment on $p$-adically
polycyclic groups in \cite[p.~347]{MR1062866}).
\end{proof}

%%%%%%%%%%%%%%%%%%%%%%%%%%%%%%%%%%%%%%%%%%%%%%%%%
\begin{Thm}
\label{thm unificador}
%%%%%%%%%%%%%%%%%%%%%%%%%%%%%%%%%%%%%%%%%%%%%%%%%
Let $X$ be a $p$-complete space such that $H^*(X; \F_p)$ is of
finite type. Assume that the module of indecomposable elements
$QH^*(X;\F_p)$ is locally finite. Then one of the following
properties is satisfied:

\begin{itemize}
\item[(1)] $X$ is aspherical,
\item[(2)] $X\langle 1 \rangle$ is equivalent to a finite product
of copies of $K(\Z^\wedge_p, 2)$,
\item[(3)] $X$ has infinitely many homotopy groups with $p$-torsion.

\end{itemize}
\noindent In the last case the space $\Omega X$ has infinitely
many non-trivial $k$-invariants.
\end{Thm}

\begin{proof}
Let $Y$ be the universal cover of $X$. Let us assume that $X$ is
not aspherical, and consider the $2$-connected cover $Y\langle 2
\rangle$ of $Y$, which can be seen as the total space in a
fibration
$$
K(\pi_2 Y, 1) \rightarrow Y\langle 2 \rangle \rightarrow Y
$$
We know that $\Omega Y$ is a $B\Z/p$-local space by
\cite[Theorem~3.2]{MR92b:55004} and $\Omega K(\pi_2 Y, 1)$ is
homotopically discrete. Therefore, if we loop once this fibration,
we see that $\Omega (Y\langle 2 \rangle)$ must be $B\Z/p$-local as
well. This means precisely that the module of indecomposable
elements $QH^* (Y\langle 2 \rangle; \F_p)$ is locally finite.
{}From Theorem~\ref{thm DW} we infer that $Y\langle 2 \rangle$
(completed at $p$) is either contractible or has infinitely many
homotopy groups with $p$-torsion, i.e. $Y$ itself has infinitely
many homotopy groups with $p$-torsion unless its $p$-completion is
an Eilenberg-Mac~Lane space of type $K(A, 2)$. In this case we
infer from Lemma~\ref{lem KZ2} that $A$ is isomorphic to a finite
direct sum of copies of~$\Z^\wedge_p$.

The statement about the $k$-invariants is a direct consequence of
Proposition~\ref{cor Lannes}. Indeed if the loop space $\Omega X$
has only a finite number of non-trivial $k$-invariants, there
exists an integer $N$ such that the $p$-complete Eilenberg-Mac
Lane space $K(\pi_n(\Omega X), n)$ is a retract of $\Omega X$ for
any $n \geq N$. Therefore this Eilenberg-Mac Lane space is
$B\Z/p$-local as well, which is only possible if $n \leq 2$. This
implies that all higher homotopy groups are trivial and so $X$ is
aspherical by the first part of the theorem.
\end{proof}

Since an unstable algebra which is locally finite as a module over
the Steenrod algebra obviously has also a locally finite module of
indecomposable elements, the Dwyer-Wilkerson condition truly
generalizes the previously accessible cases. It is in fact
straightforward to obtain the Lannes-Schwartz theorem as a
corollary.

%%%%%%%%%%%%%%%%%%%%%%%%%%%%%%%%%%%%%%%%%%%%%%%%%
\begin{Cor}
\label{cor LS}
%%%%%%%%%%%%%%%%%%%%%%%%%%%%%%%%%%%%%%%%%%%%%%%%%
Let $X$ be a simply connected space such that $H^*(X; \F_p)$ is of
finite type. Assume that $H^*(X; \F_p)$ is non-trivial and locally
finite. Then there exists an infinite number of integers $k$ such
that $\pi_k X$ contains $p$-torsion.
\end{Cor}

\begin{proof}
Since $X^\wedge_p$ is a $B\Z/p$-local space, so is its loop space
$\Omega (X^\wedge_p)$. The previous proposition applies and we can
conclude because the cohomology of $K(\Z^\wedge_p, 2)$ is not
locally finite.
\end{proof}

In view of Theorem~\ref{thm unificador} a good understanding of
the Dwyer-Wilkerson statement for arbitrary connected spaces goes
through -- compare with condition (2) -- the study of $2$-stage
Postnikov pieces. In Lemma~\ref{lem KZ2} we have identified the
second homotopy group. As for the fundamental group we will assume
that $X$ is an $\F_p$-good space, so $X^\wedge_p$ is $p$-complete.

By \cite[Proposition~3.4]{MR0438330} such spaces include all
virtually nilpotent spaces, that is, the action of the fundamental
group on any homotopy group is virtually nilpotent. Bousfield
characterizes the $H\F_p$-local spaces in
\cite[Theorem~5.5]{MR0380779} in terms of their homotopy groups,
which implies in particular that the $n$-connected covers and the
$n$-th Postnikov sections of $p$-complete (and $\F_p$-good) spaces
are $p$-complete.

In short if $X$ is a virtually nilpotent space, its $p$-completion
$X^\wedge_p$ is an $H\F_p$-local space. Its second Postnikov
section $Y = X^\wedge_p [2]$ is a $p$-complete space with only two
homotopy groups, which can be seen as the total space of a
fibration of the form
$$
K(A, 2) \rTo Y \rTo K(G, 1)
$$
where both $K(A, 2)$ and $K(G, 1)$ are $p$-complete spaces.

%%%%%%%%%%%%%%%%%%%%%%%%%%%%%%%%%%%%%%%%%%%%%%%%%
\begin{Lem}
\label{lem fundamental}
%%%%%%%%%%%%%%%%%%%%%%%%%%%%%%%%%%%%%%%%%%%%%%%%%
Let $X$ be a virtually nilpotent space such that $H^*(X; \F_p)$ is
of finite type. Then $\pi_1 (X^\wedge_p)$ is a $p$-complete group
isomorphic to $(\pi_1 X)^\wedge_p$. It is an extension of a finite
$p$-group by a nilpotent $p$-complete one.
\end{Lem}

\begin{proof}
The fundamental group $G$ of  $X^\wedge_p$ is isomorphic to that
of $K(G, 1)^\wedge_p$ by the Whitehead type theorem
\cite[Proposition~4.1]{MR1062866}. As the fundamental group of $X$
is a finitely generated virtually nilpotent group, it is in
particular polycyclic-by-finite. We conclude by Bousfield
$\F_p$-goodness result \cite[Theorem~7.2]{MR1062866} on
polycyclic-by-finite spaces that $\pi_1 (X^\wedge_p) \cong
G^\wedge_p$.

It remains to describe this virtually nilpotent group. As in
Example~\ref{ex virtuallynilpotent} we can find normal subgroups
$N \leq G' \leq G$ such that $N$ is nilpotent, finitely generated,
and torsion free, the quotient $Q=G'/N$ is $p$-perfect, and $G/G'$
is a finite $p$-group. We will actually show that the inclusion $N
\rightarrow G'$ induces an epimorphism $N^\wedge_p \rightarrow
G'^\wedge_p$. By~\cite[Lemma~5.2]{MR1062866} we only need to check
that it induces an epimorphism on the first mod $p$ homology
group, i.e. the quotient by the first term of the mod $p$ lower
central series. Notice that the quotient $N/N \cap \Gamma_1^p G'$
is isomorphic to $G'/\Gamma_1^p G'$ because $Q$ is $p$-perfect.
Therefore the maximal quotient of $N$ which is an elementary
abelian group is at least as large as $H_1(G'; \F_p)$ and we are
done.
\end{proof}

Summing up this result with Lemma~\ref{lem KZ2} we can now
describe quite accurately the $p$\guio{complete} $2$-stage Postnikov
pieces which have a $B\Z/p$-local loop space.

%%%%%%%%%%%%%%%%%%%%%%%%%%%%%%%%%%%%%%%%%%%%%%%%%
\begin{Prop}
\label{prop 2stage}
%%%%%%%%%%%%%%%%%%%%%%%%%%%%%%%%%%%%%%%%%%%%%%%%%
Let $X$ be a virtually nilpotent space such that $H^*(X; \F_p)$ is
of finite type and $\pi_n(X) = 0$ for any~$n \geq 3$. Assume that
$QH^*(X; \F_p)$ is locally finite. Then $\pi_1 (X^\wedge_p)$ is
isomorphic to $(\pi_1 X)^\wedge_p$, an extension of a finite
$p$-group by a nilpotent $p$-complete group, and $\pi_2
(X^\wedge_p)$ is isomorphic to a finite direct sum of copies
of~$\Z^\wedge_p$. \hfill{\qed}
\end{Prop}

When the fundamental group is finite, it must be a $p$-group and
we recover precisely the class of $2$-stage Postnikov systems
studied by Grodal in \cite{MR1622342}.

%%%%%%%%%%%%%%%%%%%%%%%%%%%%%%%%%%%%%%%%%%%%%%%%%
%%%%%%%%%%%%%%%%%%%%%%%%%%%%%%%%%%%%%%%%%%%%%%%%%
\section{Universal covers of $p$-completed spaces}
\label{sec groups}
%%%%%%%%%%%%%%%%%%%%%%%%%%%%%%%%%%%%%%%%%%%%%%%%%
%%%%%%%%%%%%%%%%%%%%%%%%%%%%%%%%%%%%%%%%%%%%%%%%%
In this section we consider classifying spaces in four different
families. We identify explicitly the universal covers of their
$p$-completions, and we show they are $p$-completions of spaces
inside the same family. We start by recalling the well-known case
of finite groups, see \cite{MR1308466}, even though the next
examples also contain all finite groups.

%%%%%%%%%%%%%%%%%%%%%%%%%%%%%%%%%%%%%%%%%%%%%%%%%
\subsection{Finite groups}
\label{sec finitegroups}
%%%%%%%%%%%%%%%%%%%%%%%%%%%%%%%%%%%%%%%%%%%%%%%%%
Let $G$ be a finite group and $O^p(G)$ the maximal $p$-perfect
subgroup of $G$. This is a normal subgroup and the quotient $P =
G/O^p(G)$ is a $p$-group. Therefore the fibration $B O^p(G)
\rightarrow BG \rightarrow BP$ is preserved by $p$-completion.
Since $O^p(G)$ is $p$-perfect $(BO^p(G))^\wedge_p$ is simply
connected and thus is weakly equivalent to the universal cover of
$(BG)^\wedge_p$.

Hence for any classifying space of a finite group, the universal
cover can be chosen, up to $p$-completion, to be another
classifying space.

%%%%%%%%%%%%%%%%%%%%%%%%%%%%%%%%%%%%%%%%%%%%%%%%%
\subsection{\boldmath$p$-local finite groups}
\label{sec plocalgroups}
%%%%%%%%%%%%%%%%%%%%%%%%%%%%%%%%%%%%%%%%%%%%%%%%%
Let $(S, \FF, \LL)$ be a $p$-local finite group as in
Example~\ref{ex plocal}. We learn from \cite[Theorem~4.4]{BCGLO2}
that there exists a $p$-local finite group $(O^p(S), O^p(\FF),
O^p(\LL))$ such that $|O^p(\LL)|^\wedge_p$ is the universal cover
of $|\LL|^\wedge_p$.

Again we see that the universal cover can be chosen, up to
$p$-completion, inside the class of $p$-local finite groups.

%%%%%%%%%%%%%%%%%%%%%%%%%%%%%%%%%%%%%%%%%%%%%%%%%
\subsection{Virtually nilpotent groups}
\label{sec virtnilpgroups}
%%%%%%%%%%%%%%%%%%%%%%%%%%%%%%%%%%%%%%%%%%%%%%%%%
For virtually nilpotent groups the idea to construct the universal
cover of $(BG)^\wedge_p$ out of group theoretical information is
already present in the work of Bastardas, \cite[Section
5.3]{tesisgemma}.

Let $G$ be a finitely generated virtually nilpotent group. As in
Example~\ref{ex virtuallynilpotent} we can find normal subgroups
$N \leq G' \leq G$ such that $N$ is nilpotent, finitely generated,
and torsion free, the quotient $Q=G'/N$ is $p$-perfect, and $G/G'$
is a finite $p$-group. Hence $\pi_n (BG)^\wedge_p \cong \pi_n
(BG')^\wedge_p$ for all $n \geq 2$. As we only wish to identify
the universal cover, we might as well assume from now on that $G$
sits in an extension
$$
N \rTo G \rTo Q
$$
where $Q$ is finite and $p$-perfect, and $N$ is a nilpotent group,
finitely generated, and torsion-free. In fact we see from
Lemma~\ref{lem fundamental} that $\pi_1 (BG)^\wedge_p \cong
G^\wedge_p$ is a nilpotent group since $Q$ is $p$-perfect.

The extension $N \rightarrow G \rightarrow Q$ gives rise to a
fibration $BN \rightarrow BG \rightarrow BQ$, which can be
fiberwise $p$-completed. Because $(BN)^\wedge_p$ is weakly
equivalent to the classifying space of $N^\wedge_p$, the total
space $X$ of the new fibration is the classifying space of some
group $\bar G$. Notice also that the map $BG \rightarrow B \bar G$
is a mod $p$ equivalence.

%%%%%%%%%%%%%%%%%%%%%%%%%%%%%%%%%%%%%%%%%%%%%%%%%
\begin{Lem}
\label{surjective}
%%%%%%%%%%%%%%%%%%%%%%%%%%%%%%%%%%%%%%%%%%%%%%%%%
Consider the extension $N^\wedge_p \rightarrow \bar G \rightarrow
Q$. The $p$-completion homomorphism $\bar G \rightarrow \bar
G^\wedge_p$ is then surjective.
\end{Lem}

\begin{proof}
In the proof of Lemma~\ref{lem fundamental} we showed that the
inclusion $N \hookrightarrow G$ induces an epimorphism $N^\wedge_p
\rightarrow G^\wedge_p$. Since $G^\wedge_p \cong \bar G^\wedge_p$,
the morphism $\bar G \rightarrow \bar G^\wedge_p$ is onto as well.
\end{proof}

Let us define now $K$ as the kernel  of the completion morphism
$\bar G \rightarrow \bar G^\wedge_p$ (the intersection of all the
terms in the mod $p$ lower central series). In
Theorem~\ref{universalcover} we identify the universal cover of
$(BG)^\wedge_p$  as the $p$-completed classifying space of the
group~$K$.

%%%%%%%%%%%%%%%%%%%%%%%%%%%%%%%%%%%%%%%%%%%%%%%%%
\begin{Prop}
\label{barfibration}
%%%%%%%%%%%%%%%%%%%%%%%%%%%%%%%%%%%%%%%%%%%%%%%%%
There is a fibration $(BK)^\wedge_p \rTo (B \bar G)^\wedge_p \rTo
B(G^\wedge_p)$.
\end{Prop}

\begin{proof}
Let us consider the pair of fibrations  $B\bar G \rightarrow BQ$
and $B\bar G \rightarrow B(G^\wedge_p)$ with fibers
$(BN)^\wedge_p$ and $BK$ respectively. We deduce from
Lemma~\ref{surjective} that the induced map on fundamental groups
$\bar G \rightarrow G^\wedge_p \times Q$ is onto. Since the first
fiber $(BN)^\wedge_p$ is nilpotent, we can apply the ``nilpotent
action lemma" \cite[5.1]{MR0438330} of Dwyer, Farjoun, and Kan.
Therefore the nilpotent group $G^\wedge_p$ acts nilpotently on all
homology groups $H_n(BK; \F_p)$. The ``nilpotent fibration lemma"
\cite[Proposition~4.2(i)]{MR0438330} tells us now that
$p$-completion preserves the fibration $BK \rightarrow B\bar G
\rightarrow B(G^\wedge_p)$ and we are done.
\end{proof}

%%%%%%%%%%%%%%%%%%%%%%%%%%%%%%%%%%%%%%%%%%%%%%%%%
\begin{Thm}
\label{universalcover}
%%%%%%%%%%%%%%%%%%%%%%%%%%%%%%%%%%%%%%%%%%%%%%%%%
Let $G$ be a virtually nilpotent group which is an extension of a
finite, $p$\guio{perfect} group by a nilpotent, finitely generated, and
torsion-free one. The space $(BK)^\wedge_p$ is the universal cover
of $(BG)^\wedge_p$.
\end{Thm}

\begin{proof}
One knows that the fundamental group of $(BG)^\wedge_p$ is
$G^\wedge_p$ (\cite[Theorem~7.2]{MR1062866}). In view of the above
proposition we only need to remark that $BG \rightarrow B \bar G$
is an $H\F_p$-equivalence (the $p$-completion $N \rightarrow
N^\wedge_p$ is so for $N$ is nilpotent).
\end{proof}

We note that $K$ is virtually nilpotent as well, being a subgroup
of a virtually nilpotent one. We can even say more, since
$N^\wedge_p$ is the $p$-completion of a finitely generated torsion
free nilpotent group: $K$ is $p$-adically polycyclic-by-finite. In
the situation where $G$ is actually finite, $K$ coincides with
$O^p(G)$, which is consistent with the approach of
Levi~\cite{MR1308466}, compare with~\ref{sec finitegroups}.

%%%%%%%%%%%%%%%%%%%%%%%%%%%%%%%%%%%%%%%%%%%%%%%%%
\subsection{Quasi \boldmath$p$-perfect groups}
\label{sec quasiperfectgroups}
%%%%%%%%%%%%%%%%%%%%%%%%%%%%%%%%%%%%%%%%%%%%%%%%%
By definition, see Example~\ref{ex quasiperfect}, a group $G$ is
quasi $p$-perfect if the subgroup $\Gamma_1^p(G)$ is $p$-perfect.
The fibration
$$
B\Gamma_1^p(G) \rTo BG \rTo BH_1(G; \F_p)
$$
is preserved under $p$-completion since $H_1(G; \F_p)$ is an
(elementary abelian) $p$-group and $BG$ is $\F_p$-good. Hence the
universal cover of $(BG)^\wedge_p$ is $\bigl( B\Gamma_1^p(G)
\bigl)^\wedge_p$.

%%%%%%%%%%%%%%%%%%%%%%%%%%%%%%%%%%%%%%%%%%%%%%%%%
%%%%%%%%%%%%%%%%%%%%%%%%%%%%%%%%%%%%%%%%%%%%%%%%%
\section{Homotopy groups of $p$-completed classifying spaces}
\label{sec homotopygroups}
%%%%%%%%%%%%%%%%%%%%%%%%%%%%%%%%%%%%%%%%%%%%%%%%%
%%%%%%%%%%%%%%%%%%%%%%%%%%%%%%%%%%%%%%%%%%%%%%%%%
We obtain in this section in a single proof the results which were
known before about $p$\guio{completions} of classifying spaces of finite
groups (Levi) and virtually nilpotent groups (Bastardas). We prove
along the same lines a new result for quasi $p$-perfect groups and
$p$-local finite groups.

%%%%%%%%%%%%%%%%%%%%%%%%%%%%%%%%%%%%%%%%%%%%%%%%%
\begin{Thm}
\label{groupdichotomy}
%%%%%%%%%%%%%%%%%%%%%%%%%%%%%%%%%%%%%%%%%%%%%%%%%
Let $X$ be the classifying space of a member of the
 following four families:
\begin{itemize}
\item[(1)]  finite groups,

\item[(2)]  $p$-local finite groups,

\item[(3)]  finitely generated virtually nilpotent groups,

\item[(4)]  quasi $p$-perfect groups of finite virtual mod $p$
cohomology.

\end{itemize}

Then the $p$-completion of $X$ is either aspherical or it has
infinitely many homotopy groups with $p$-torsion. In this case the
space $\Omega (X)^\wedge_p$ has infinitely many non-trivial
$k$\guio{invariants}.

\end{Thm}

\begin{proof}
In view of the previous section it remains to prove that none of
the four families can contain a space whose $p$-completion is
$K(\Z^\wedge_p, 2)$. The mod $p$ cohomology of such a space is
polynomial on a generator in dimension~$2$, hence concentrated in
even dimensions. Therefore there are no Bocksteins at all, which
means by an elementary Bockstein spectral sequence argument that
there is no $p$-torsion in the higher integral homology groups.

When $G$ is a finite group, choose $L$ to be the trivial subgroup.
When $G$ is a virtually nilpotent group, we can assume as
in~\ref{sec virtnilpgroups} that $G=K$ is $p$-perfect and so
$(BG)^\wedge_p$ is simply connected. Choose then $L$ to be the
nilpotent subgroup of finite index $K \cap N^\wedge_p$. This is a
subgroup of $N^\wedge_p$, which has finite homological dimension,
\cite[VIII.2]{MR672956}. Finally when $G$ is a quasi $p$-perfect
group of virtually finite mod $p$ cohomology, we choose $L$ to be
some subgroup of finite index which has finite mod $p$ cohomology.
A standard transfer argument applied to the subgroup $L < G$ shows
now that in all three cases the multiplication by the (finite)
index is zero on high enough integral homology groups of $G$.

In case (2) the mod $p$ cohomology of a $p$-local finite group is
contained as a retract in the cohomology of its Sylow $p$-subgroup
as a unstable subalgebra over the Steenrod algebra
\cite[Proposition~5.5]{MR1992826} (see also
\cite[Proposition~9.4]{Kari}). The cohomology of $K(\Z^\wedge_p,
2)$ which has no (higher) Bocksteins cannot be a retract of the
cohomology of a finite group.
\end{proof}

We point out that the assumptions made on the virtually nilpotent
group could hardly be relaxed. For example, if one drops the
finitely generated hypothesis, the result is obviously false, as
shown by the example of the Pr\"ufer group:
$(B\Z/p^\infty)^\wedge_p \simeq K(\Z^\wedge_p,2)$.

Let us mention that the $p$-completed classifying spaces of
$p$-perfect groups of finite virtual cohomological dimension have
been studied by Levi. He proves in \cite[Theorem~1.4]{MR1709097}
that $\Omega (BG)^\wedge_p$ is a retract of some finite complex.

%%%%%%%%%%%%%%%%%%%%%%%%%%%%%%%%%%%%%%%%%%%%%%%%%%
%%%%%%%%%%%%%%%%%%% REFERENCES %%%%%%%%%%%%%%%%%%%
%%%%%%%%%%%%%%%%%%%%%%%%%%%%%%%%%%%%%%%%%%%%%%%%%%

\providecommand{\bysame}{\leavevmode\hbox
to3em{\hrulefill}\thinspace}
\providecommand{\MR}{\relax\ifhmode\unskip\space\fi MR }
\providecommand{\MRhref}[2]{%
  \href{http://www.ams.org/mathscinet-getitem?mr=#1}{#2}
} \providecommand{\href}[2]{#2}

%%%%%%%%%%%%%%%%%%%%%%%%%%%%%%%%%%%%%%

%%%%%%%%%%%%%%%%%%%%%%%%%%%%%%%%%%%%%%%%%%%%%%%%%
%%%%%%%%%%%%%%%%%%% ADDRESSES %%%%%%%%%%%%%%%%%%%
%%%%%%%%%%%%%%%%%%%%%%%%%%%%%%%%%%%%%%%%%%%%%%%%%

\bigskip
{\small
\begin{center}
\begin{minipage}[t]{8 cm}
Nat\`{a}lia Castellana and J\'er\^{o}me Scherer\\ Departament de Matem\`atiques,\\
Universitat Aut\`onoma de Barcelona,\\ E-08193 Bellaterra, Spain
\\\textit{E-mail:}\texttt{\,natalia@mat.uab.es}, \\
\phantom{\textit{E-mail:}}\texttt{\,jscherer@mat.uab.es}
\end{minipage}
\begin{minipage}[t]{7 cm}
Juan A. Crespo \\ Departament de Economia i de Hist\`{o}ria
Econ\`{o}mica,
\\ Universitat Aut\`onoma de Barcelona,\\ E-08193 Bellaterra,
Spain
\\\textit{E-mail:}\texttt{\,juanalfonso.crespo@uab.es}
\end{minipage}
\end{center}}

%%%%%%%%%%%%%%%%%%%%%%%%%%%%%%%%%%%%%%
\end{document}